\documentclass[12pt,a4paper]{article}

\usepackage{type1cm}        % activate if the above 3 fonts are
\usepackage[left=20mm,right=20mm, top=23mm, bottom=23mm]{geometry}
\usepackage{makeidx}         % allows index generation
\usepackage{graphicx}        % standard LaTeX graphics tool
\usepackage{multicol}        % used for the two-column index
\usepackage[bottom]{footmisc}% places footnotes at page bottom

\usepackage{amsmath}
\usepackage{amsfonts}
\usepackage{amssymb}
\usepackage{mathtools} % must be loaded
\usepackage{amsthm}
\RequirePackage{tikz}
\usetikzlibrary{arrows.meta}
\usepackage{pgfplots} 
\usepackage{pgfgantt}
\usepackage{pdflscape}
\usepackage{bm}
\pgfplotsset{compat=newest} 
\pgfplotsset{plot coordinates/math parser=false}

\usepackage{newtxtext}       % 
\usepackage[varvw]{newtxmath}       % selects Times Roman as basic font

\makeindex             % used for the subject index
\newcommand{\Dx}{\Delta x}
\newcommand{\Dt}{\Delta t}
\newcommand{\N}{\mathbb{N}}
\newcommand{\R}{\mathbb{R}}
\newcommand{\Z}{\mathbb{Z}}
\newcommand\half{{\frac12}}

\newcommand\jmh{{j-\half}}

\newcommand\jph{{j+\half}}
\newcommand\jpo{{j+1}}

\newcommand{\mm}{\operatorname{minmod}}
\newcommand{\dr}{\rho^{\Dx}}
\newcommand{\BV}{\operatorname{BV}}
\newcommand{\ndt}{{\eta}}
\newcommand{\wt}{\omega_{\ndt}}

\newcommand{\Ne}{N_\ndt}
\newcommand{\norm}[1]{{\left\|#1\right\|}}
\newcommand{\abs}[1]{{\left|#1\right|}}
\newcommand\brho{{\bm{\rho}}}
\newlength\fwidth

\theoremstyle{plain}						% Default
\newtheorem{theorem}{Theorem}[section]

\newtheorem{corollary}[theorem]{Corollary}
\theoremstyle{definition}

\theoremstyle{remark}
\newtheorem{remark}[theorem]{Remark}

\newcommand{\rv}[1]{{#1}}
\begin{document}

\title{A note on the central-upwind scheme for nonlocal conservation laws}
\author{Jan Friedrich\footnotemark[1], Samala Rathan\footnotemark[2] and Sanjibanee Sudha\footnotemark[2]}

\footnotetext[2]{Chair of Optimal Control, Department for Mathematics, School of Computation, Information and Technology, Technical University of Munich, Boltzmannstraße 3, 85748 Garching b. Munich, Germany, jan.friedrich@cit.tum.de}
\footnotetext[1]{Department of Humanities and Sciences, Indian Institute of Petroleum and Energy, Visakhapatnam, Andhra Pradesh, India-530003 (\{rathans.math,sudhamath21\}@iipe.ac.in)}
\maketitle

\begin{abstract}
    The central-upwind flux is a widely used numerical flux function for local conservation laws. It has been investigated by Kurganov and Polizzi (2009) for a specific nonlocal conservation law and can be derived from a fully-discrete second-order scheme. Here, we derive this fully-discrete scheme in detail with a particular focus on the occurring nonlocal terms. In addition, we derive the central-upwind flux for a class of nonlocal conservation laws and use an estimate on the nonlocal speed which fixes the nonlocality at the cell interfaces. We prove that the resulting first-order numerical scheme converges to the correct solution. Under additional assumptions on the analytical flux we present a similar result for a second-order central-upwind scheme. Numerical examples compare the central-upwind schemes to Godunov-type schemes and the fully-discrete scheme.
\end{abstract}

\section{Introduction}
\label{sec:1}
In this work we focus on one-dimensional scalar conservation laws with a nonlocality in space.
Such models can be found in various physical applications such as supply chains~\cite{goettlich2010supplychains}, sedimentation~\cite{betancourt2011nonlocal}, and traffic flow~\cite{BlandinGoatin2016, chiarello2018global, friedrich2020nonlocal, friedrich2018godunov}.
They are described by
\begin{align}\label{eq:conslaw}
\partial_t \rho(t,x)+\partial_x(F(\rho,\wt \ast \rho))(t,x)&=0,\quad &(t,x)\in\mathbb{R}_+\times\mathbb{R}, \nonumber\\
    \rho(x,0) &=\rho_0(x), \quad &x\in\mathbb{R},  
\end{align}
where $\rho$ is the state variable, $t$ the time, $x$ the space variable, $F$ is a suitable flux function and the convolution term  $\wt\ast\rho$ is an integral evaluation over space.
Here, we consider a forward looking kernel with compact support on $[0,\eta]$:
\begin{align*}
    (\wt\ast \rho) (t,x):=\int_x^{x+\eta} \wt(y-x)\rho(t,y) dy.
\end{align*}
A well studied flux $F$ is given by
\begin{equation}\label{eq:Fsimple}
  F(\rho,\wt \ast \rho)=g(\rho)v(\wt \ast \rho).
\end{equation}
Under the following hypotheses the well-posedness of \eqref{eq:conslaw}--\eqref{eq:Fsimple}, i.e. existence and uniqueness of weak entropy solutions, has been shown in \cite{chiarello2018global}
\begin{enumerate}
    \item $\rho(0,x)=\rho_0(x) \in \BV(\mathbb{R},I), \, \, I={[\rho_{\min},\rho_{\max}]}\subseteq\mathbb{R^+},$
    \item $g\in C^1(I;\mathbb{R}^+),$
    \item $v\in C^2(I;\mathbb{R}^+),\quad v'\leq 0,$
    \item  $\omega_{\eta}\in C^1([0,\eta]),\mathbb{R^+})$\,\,\,\, \text{with}\,\,\, $\wt^{\prime}\leq0,$\,\,\,$\int_{0}^{\eta}\omega_{\eta}(x)dx=1$ \,\,\,$\forall\eta>0,$\,\,\,\,$\underset{\eta\to\infty}{\lim} \omega_{\eta}(0)=0,$
\end{enumerate}
where $\rho_{\min}:=\underset{\mathbb{R}}{\min}(\rho_0)$ and $\rho_{\max}:=\underset{\mathbb{R}}{\max}(\rho_0)$ \rv{and $BV(\R;I)$ denotes the set of all functions with a bounded total variation mapping from $\R$ to $I$.}
We note that for the flux \eqref{eq:Fsimple} weak entropy solutions are necessary however if $g$ is linear, weak solutions are already unique \cite{KeimerPflug2017}.

In contrast to local conservation laws, no Riemann solvers can be constructed due to the global dependence of the flux on the solution.
This makes it difficult to derive Godunov-type numerical schemes, however they still can be constructed under certain assumptions \cite{friedrich2018godunov,friedrich2023numerical}.
\rv{If $g$ is linear,} Godunov-type schemes coincide with Upwind-type schemes and provide a solid and accurate approximation of \eqref{eq:conslaw}--\eqref{eq:Fsimple}.
However, once $g$ becomes nonlinear central methods might be preferable, as e.g. the Lax-Friedrichs-type flux.
Central schemes can be used as black-box solvers for hyperbolic (local and nonlocal systems of) conservation laws.
However, they usually suffer from excessive numerical dissipation.
For nonlocal conservation laws such finite volume numerical methods (of first- and second-order) based on the Lax-Friedrichs and Nessyahu-Tadmor scheme have been studied in \cite{GoatinScialanga2016,kurganov2009non,sudha2025convergence}.
One way to reduce the numerical diffusion is to use higher-order methods \cite{ChalonsGoatinVillada2018,friedrich2019maximum}.

Another approach is to derive more accurate central numerical flux functions.
This approach is very popular for local conservation laws and was first introduced in \cite{kurganov2000new}. Here the central Kurganov-Tadmor (KT) scheme was developed to derive the so called central-upwind (CU) scheme.
Further, the approach was extended in \cite{kurganov2007reduction}.
The KT scheme is a second-order accurate method that incorporates local maximum wave speeds (related to the CFL condition) at cell boundaries. This significantly reduces the numerical viscosity. Moreover, expressing it in a semi-discrete form, leads to the CU flux which inherits the central nature, includes a build-in anti-diffusion term and allows for higher-order extensions.
The CU flux has been considered by many authors for a broad range of hyperbolic partial differential equations.
We refer to \cite{kurganov2000new,kurganov2007reduction,kurganov2009non,cui2024bound,kurganov2001centralupwind} and the references therein.\\
Back in 2009, the CU flux for nonlocal conservation laws was already derived via the KT scheme in \cite{kurganov2009non}.
However, the actual fully-discrete scheme was only used as an intermediate step without presenting its details, in particular concerning the nonlocality, and both schemes have only been considered for a specific choice of a nonlocal flux function  (i.e. a specific choice of $g$ and $v$) and constant or linear kernels. In addition, no theoretical results were presented.
Here, we revisit the work \cite{kurganov2009non} and extend the approaches to the more general flux \eqref{eq:Fsimple} and nonincreasing kernel functions.
Further, we discuss the convergence and properties of first- and second-order CU schemes.
One main difference to \cite{kurganov2009non} is a more precise estimate of the nonlocal speeds, which is detailed below.

Thereby, the work is structured as following:
In Section \ref{sec:2}, we derive the fully-discrete scheme with a emphasize on the nonlocal terms. Then, we use it to obtain the CU flux and prove the convergence of the resulting first-order scheme for the nonlocal conservation law \eqref{eq:conslaw}--\eqref{eq:Fsimple}. In addition, we present a result for the second-order CU scheme, Further, we discuss how to extend the schemes to systems of nonlocal balance laws. In Section \ref{sec:3}, we compare the schemes numerically. Finally, we provide conclusions in Section \ref{sec:4}.

\section{The Kurganov-Tadmor and central-upwind scheme for nonlocal conservation laws}
\label{sec:2}
As aforementioned, we start by deriving a second-order fully-discrete scheme.
Therefore we integrate the conservation law \eqref{eq:conslaw} over a smooth and non-smooth area. This results in the KT scheme which can be used to obtain the CU flux.
The approach to derive both schemes follows very closely the ones presented in \cite{kurganov2007reduction,kurganov2001centralupwind,kurganov2000new}.
The key is to estimate the maximal and minimal speed at the cell interface to obtain smooth and non-smooth areas.
Apart from extending the approach from \cite{kurganov2009non} to a more general class of conservation laws, the main difference is the estimate of the nonlocal speeds, i.e. instead of an upper bound on the nonlocal term, we use its approximation at the cell interface.\\
To be more precise we discretize \eqref{eq:conslaw} in space and time by an equidistant grid, where $\Dx$ is the step size in space and $\Dt$ is the step size in time.
Hence, $t^n=n\Dt$ with $n\in\N$ describes the time grid and $x_j=j\Dx,\ j\in\Z$ the cell centers of the space grid with the cell interfaces $x_{j\pm 1/2}$.
Let, $\rho_{j+1/2}^\pm$ be the left and right approximations of the approximate solution at $x_{j+1/2}$ and time $t^n$, e.g. for a piece-wise constant approximation, we have $\rho_{j+1/2}^-=\rho_j$ and $\rho_{j+1/2}^+=\rho_{j+1}$. 
At the cell interface the nonlocal term is, as an integral evaluation, continuous in $x$, such that no left or right limit is needed.
In particular, we have
\begin{align*}
    &(\wt \ast \rho)(t,x_\jph)= \int_{x_\jph}^{x_\jph+\ndt} \wt(y-x_\jph)\rho(t,y)dy
\end{align*}
We assume that for fixed $R\in\mathbb{R}$ the flux $F(\rho,R)$ is either a convex or concave function in $\rho$. Then, the estimated maximum and minimum speeds at $x_\jph$ are given by
$$c_{j+1/2}^{+}=\max\left(\frac{\partial F(\rho_{j+1/2}^-,(\wt \ast \rho)(t,x_\jph))}{\partial \rho},\frac{\partial F(\rho_{j+1/2}^+,(\wt \ast \rho)(t,x_\jph))}{\partial \rho},0\right)$$
and
$$c_{j+1/2}^{-}=\min\left(\frac{\partial F(\rho_{j+1/2}^-,(\wt \ast \rho)(t,x_\jph))}{\partial \rho},\frac{\partial F(\rho_{j+1/2}^+,(\wt \ast \rho)(t,x_\jph))}{\partial \rho},0\right)$$
We note that in the following the derivation of both schemes does not require the flux function to have the form \eqref{eq:Fsimple}, however for the convergence of the CU scheme we only study \eqref{eq:Fsimple}.

\subsection{The Kurganov-Tadmor scheme}
\begin{figure}[h]
\centering
% \documentclass{standalone}
% \usepackage{tikz}
% \usepackage{amsmath}
% Adjust sign of a number for the tikz picture of the grid 
	\newcommand\checkindexnew[1]{
		\pgfmathsetmacro{\var}{#1}
		\pgfmathparse{ifthenelse(\var==0, "",ifthenelse(\var>0, "+#1","#1"))} \pgfmathresult}%

% \begin{document}
\begin{tikzpicture}[scale=1]
% Draw the x-axis
    % \draw[->] (-5, 0) -- (4, 0) node[right] {$x$}; % x-axis from -4 to 4, with an arrow
    
    % % Draw the y-axis
    % \draw[->] (-5, 0) -- (-5, 4) node[above] {$y$}; % y-axis from -4 to 4, with an arrow

  % Draw x-axis with an arrow
  \draw[->] (-5.5, -5) -- (3.5, -5) node[right] {$x$};

  % Define coordinates for grid points
  \coordinate (xjm32r) ;    % x_{j-3/2,r}
  \coordinate (xjm1); %at (-4.3, -5);    % x_{j-1}
  \coordinate (xjm12l); %at (-3.6, -5);  % x_{j-1/2,l}
  \coordinate (xjm12); %at (-2.7, -5);   % x_{j-1/2}
  \coordinate (xjm12r);%at (-1.8, -5);  % x_{j-1/2,r}
  \coordinate (xj); %at (-1.1, -5);      % x_j
  \coordinate (xjp12l); %at (-0.4, -5);  % x_{j+1/2,l}
  \coordinate (xjp12); %at (0.5, -5);    % x_{j+1/2}
  \coordinate (xjp12r); %at (1.6, -5);   % x_{j+1/2,r}
  \coordinate (xjp1);% at (2.3, -5);     % x_{j+1}
  \coordinate (xjp32l); %at (3, -5);     % x_{j+3/2,l}
\coordinate (xjp32); %at (3, -5);     % x_{j+3/2,l}
  
%  \draw (xjm32r) (-5, -5) node[below] {$x_{j-\frac{3}{2},r}^n$};  % Label
\draw (xjm1) (-4.3, -5) node[below] {$x_{j-1}$};  % Label
\draw (xjm12l) (-3.6, -5) node[below] {$x_{j-\frac{1}{2},l}^n$};  % Label
  \draw (xjm12) (-2.7, -5) node[below] {$x_{j-\frac{1}{2}}$};  %(xjm12) at % Label
  \draw (xjm12r) (-1.8, -5) node[below] {$x_{j-\frac{1}{2},r}^n$};

  \draw (xj) (-1.1, -5) node[below] {$x_{j}$};
  \draw (xjp12l) (-0.4, -5) node[below] {$x_{j+\frac{1}{2},l}^n$};  % Label
\draw (xjp12) (0.5, -5) node[below] {$x_{j+\frac{1}{2}}$};  % Label
\draw (xjp12r) (1.6, -5) node[below] {$x_{j+\frac{1}{2},r}^n$};  % Label
 \draw (xjp1) (2.3, -5) node[below] {$x_{j+1}$};  % Label
%  \draw (xjp32l) (3, -5) node[below] {$x_{j+\frac{3}{2},l}^n$};  % Label
%\draw (xjp32) (3.5, -5) node[below] {$x_{j+\frac{3}{2}}$};  % Label

% Adding labels for the functions at points
\node at (-4.1, -3.3) {${\rho}^{n}_{j-1}$};
\node at (-1.1, -4) {${\rho}^{n}_{j}$};
\node at (2, -3.1) {${\rho}^{n}_{j+1}$};
% Solid sloped lines
\draw[thick] (-5.5, -3) -- (-2.7, -3);
\draw[thick] (-2.7,-3.5) -- (0.5,-3.5);  
\draw[thick] (0.5,-2.8) -- (3.5,-2.8);
% Draw an inclined line through the midpoint with slope 1 (45 degrees)
    \draw[thick] (-5.5, -2.5) -- (-2.7, -3.4); % Adjust endpoints if needed
 % \draw[thick] (-2.7, -3.2) -- (0.5, -3.7); % Adjust endpoints if needed
\draw[thick] (0.5,-3) -- (3.5, -2.5);
%Dashed line 
\draw[dotted] (-2.7,-1.8)--(-2.7,-5);
\draw[dotted] (0.5,-1.6)--(0.5, -5);
%Draw a thick line second level (smooth)
\draw[thick] (-5.5,-1.4)--(-3.6,-1.4);
\draw[thick] (-1.8,-2.2)--(-0.4,-2.2);
\draw[thick] (1.6,-1)--(3.5,-1);
%Draw a thick line second level (Non-smooth)
\draw[thick] (-3.6,-1.8)--(-1.8,-1.8);
\draw[thick] (-0.4,-1.6)--(1.6,-1.6);
%Draw a dashed line smooth line
\draw[dashed] (-0.4,-5)--(-0.4,-1.6);
\draw[dashed] (1.6,-5)--(1.6,-1);
\draw[dashed] (-3.6,-5 )--(-3.6,-1.4);
\draw[dashed] (-1.8,-5)--(-1.8,-1.8);
%\draw[dashed] (1.6,-0.6)--(1.6,-3);

% Adding labels for the functions at points
\node[above] at (-4.3, -1.4) {$w^{n+1}_{j-1}$};
\node[above] at (-1.1, -2.2) {$w^{n+1}_{j}$};
\node[above] at (2.3, -1) {$w^{n+1}_{j+1}$};
%midpoint smooth part order slope line second level
 \draw[thick] (-3.6, -1.6) -- (-1.8, -2); % Adjust endpoints if needed
 \draw[thick] (-0.4, -1.8) -- (1.6, -1.4); % 
% Adding labels for the functions at points(smooth level)
\node at (-2.7, -2.1) {${w}^{n+1}_{j-1/2}$};
\node at (0.5, -1.9) {${w}^{n+1}_{j+1/2}$};
%Final step
\draw[thick](-2.7,-0.5)--(0.5,-0.5);
\draw[thick](-2.7,-0.5)--(-2.7,-1.8);
\draw[thick](0.5,-0.5)--(0.5,-1.6);
% Adding labels for the functions at points(final step)
\node[below] at (-1.1, -0.5) {$\rho^{n+1}_{j}$};

% \documentclass{standalone}
% \usepackage{tikz}
% \usepackage{amsmath}

% \begin{document}

% \begin{tikzpicture}[scale=1.3, >=stealth]

% % --- Axes ---
% \draw[->] (-5.5,-4.5) -- (3.5,-4.5) node[right] {$x$};

% % --- Cell interfaces and centers ---
% \foreach \x/\lbl in {-2.7/{x_{j-\frac{1}{2}}},0.5/{x_{j+\frac{1}{2}}}} {
%   \draw[dashed] (\x,-4.5) -- (\x,-1);
%   \node[below] at (\x,-4.5) {$\lbl$};
% }
% \node[below] at (-4.3,-4.5) {$x_{j-1}$};
% \node[below] at (-1.1,-4.5) {$x_{j}$};
% \node[below] at (2.3,-4.5) {$x_{j+1}$};

% % --- Piecewise constant reconstruction (flat segments) ---
% \draw[thick] (-5.5,-3.0) -- (-2.7,-3.0);
% \draw[thick] (-2.7,-3.5) -- (0.5,-3.5);
% \draw[thick] (0.5,-2.8) -- (3.5,-2.8);

% % --- Vertical jumps between cells ---
% \draw[thick] (-2.7,-3.0) -- (-2.7,-3.5);
% \draw[thick] (0.5,-3.5) -- (0.5,-2.8);

% % --- Cell-average labels ---
% \node at (-4.1,-2.8) {$\rho_{j-1}^{n}$};
% \node at (-1.1,-3.3) {$\rho_{j}^{n}$};
% \node at (2.0,-2.6) {$\rho_{j+1}^{n}$};

% % --- Next time step (piecewise constant update) ---
% \draw[thick,blue] (-5.5,-1.5) -- (-2.7,-1.5);
% \draw[thick,blue] (-2.7,-1.0) -- (0.5,-1.0);
% \draw[thick,blue] (0.5,-1.3) -- (3.5,-1.3);
% \draw[thick,blue] (-2.7,-1.5)--(-2.7,-1.0);
% \draw[thick,blue] (0.5,-1.0)--(0.5,-1.3);

% % --- Updated labels ---
% \node[blue] at (-4.1,-1.3) {$\rho_{j-1}^{n+1}$};
% \node[blue] at (-1.1,-0.9) {$\rho_{j}^{n+1}$};
% \node[blue] at (2.0,-1.2) {$\rho_{j+1}^{n+1}$};

% % --- Legend or label ---
% \node[below right, align=left] at (-5.5,-2.2) {\small piecewise constant reconstruction};

\end{tikzpicture}

 % \end{document}
\caption{Grid and reconstructions used during the fully-discrete  scheme.}
\label{fig:tikz}
\end{figure}
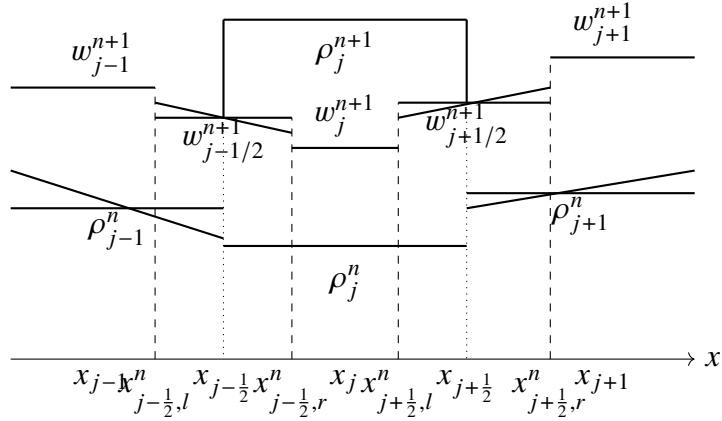

Using the estimates of the speed at each cell interface we can split cells in smooth and non-smooth parts.
In particular, the points $x_{j+1/2,l}^n=x_{j+1/2}+c_{j+1/2}^{-}\Delta t$ and $x_{j+1/2,r}^n=x_{j+1/2}+c_{j+1/2}^{+}\Delta t$ separate between these regions, see Figure \ref{fig:tikz}. Note that the speed estimates $c_\jph^\pm$ depend on time, too. We omit the super index here, for readability.
The main idea is to integrate over the non-smooth region $[x_{j+1/2,l}^n,x_{j+1/2,r}^n]$ and smooth region $[x_{j-1/2,r}^n,x_{j+1/2,l}^n]$ to derive a second-order accurate approximation.
Finally, these solutions are projected back to the original grid.
Thereby, the derivation is rather similar to the local case discussed in \cite{kurganov2007reduction} and can be done under the CFL condition
\begin{align}\label{eq:CFL0}
  \Delta t \cdot \max_j \left( \max \left( c^+_{j+\frac{1}{2}}, -c^-_{j+\frac{1}{2}} \right) \right) \leq \frac{\Delta x}{2},
\end{align}
\rv{which guarantees that the smooth and non-smooth regions do not overlap.}
We only show the resulting approximations and discuss the approximations of the different nonlocal terms in detail.
We assume that the cell-averages of the solution are given at a time $t^n$ by
\begin{equation}\label{eq:cellavg}
  \rho_{j}^n=\frac{1}{\Delta x}\int_{x_{j-1/2}}^{x_{j+1/2}}\rho(t^n,x)dx.  
\end{equation}
To obtain a second-order accurate approximation we reconstruct a piecewise linear function 
\begin{equation}\label{eq:piecewise}
    \rho^{\Delta x}(t,x)={\rho}_{j}^n+{s_{j}^n}(x-x_j) \,\,\text{for}\,\, (t,x)\in[t^n,t^{n+1})\times [x_{j-1/2},x_{j+1/2}).
\end{equation}
Here, the slopes $s_j^n$ approximate the spatial derivative of $\rho$ and are computed via the minmod limiter
\begin{align}\label{eq:slope for X}
s_{j}&=\mm\bigg(\frac{\rho_{j}^n-\rho_{j-1}^n}{\Delta x},\frac{\rho_{j+1}^n-\rho_{j}^n}{\Delta x}\bigg),\\
\quad\mm(a, b) &= 
\begin{cases} 
      a, & \text{if } |a| < |b| \text{ and } a \cdot b > 0, \\
      b, & \text{if } |b| < |a| \text{ and } a \cdot b > 0, \\
      0, & \text{if } a \cdot b \leq 0. 
\end{cases}\nonumber
\end{align}
Hence, we obtain the following left and right values at $x_\jph$ used to compute $c_\jph^\pm$
\begin{equation*}
\rho_{j+1/2}^{+}=\rho_{j+1}^{n}-\frac{\Delta x}{2}s_{j+1}^{n},\quad \rho_{j+1/2}^{-}=\rho_j^{n}+\frac{\Delta x}{2} s_j^{n}.
\end{equation*}
To estimate the nonlocal speeds, it remains to approximate the nonlocal term at $x_{j+1/2}$. Here, we choose the following second-order accurate approximation:
\begin{equation}\label{eq:approximatedNLterm}
\begin{aligned}
    &(\wt \ast \rho)(t^n,x_\jph)= \int_{x_\jph}^{x_\jph+\ndt} \wt(y-x)\rho(t^n,y)dy\\
    &\approx \sum_{k=0}^{\Ne-1} \gamma_k \rho_{j+k+1}^n+\gamma_{\Ne+1}(\rho_{j+\Ne+1}^n-s_{j+\Ne+1}^n(x_{\jph+\ndt}-x_{j+\Ne+1}))=:R_\jph^n\\
&\text{with}\quad 
    \gamma_k=\int_{k\Dx}^{\min\{(k+1)\Dx,\eta\}} \wt (x) dx\quad \text{and}\quad \Ne:=\lfloor \ndt/\Dx\rfloor.
\end{aligned}
\end{equation}
For brevity, we use the notation $R_\jph^n\approx (\wt\ast\rho)(t^n,x_\jph)$ in the following.
\begin{remark}
    We note that for the accuracy of the numerical scheme, it is not necessary to use a second-order accurate approximation of the nonlocal term, since we are only considering an estimate of the speeds.
    However, better estimates result in more accurate schemes, e.g. it is also possible to compute the nonlocal term over the piecewise linear function exactly.
    The above quadrature rule (instead of the midpoint rule) is chosen, since it is consistent with the first-order scheme detailed below.
\end{remark}
The integration of the piecewise linear function over the non-smooth area
$[x_{j+1/2,l}^n,$ $x_{j+1/2,r}^n]\times[t^n,t^{n+1}]$ gives the following intermediate cell-averages at $t^{n+1}$ 
\begin{equation*}
    \begin{aligned}
        w_{j+1/2}^{n+1} =&\frac{1}{c_{j+1/2}^{+}-c_{j+1/2}^{-}}\bigg[\rho_{\jph,r}^nc_{j+1/2}^{+}-\frac{s_\jpo^n}{2}(c_\jph^+)^2\Delta t-\rho_{\jph,l}^nc_{j+1/2}^{-}+\frac{s_j^n}{2}(c_\jph^-)^2\Delta t\\
       -&(F(\rho_{j+1/2,r}^{n+1/2},R_{j+1/2,r}^{n+1/2})-F(\rho_{j+1/2,l}^{n+1/2},R_{j+1/2,l}^{n+1/2}))\bigg].
    \end{aligned}
\end{equation*}
Similarly, for the smooth area we obtain
\begin{equation*}
    \begin{aligned}
       w_j^{n+1}=&{\rho}_j^n+\frac{s_j^n}{2}(c_\jph^++c_\jph^-)
       -\frac{\Delta t}{\Delta x-\Delta t(c_{j-1/2}^{+}-c_{j+1/2}^{-})}\bigg[F(\rho_{j+1/2,l}^{n+1/2},R_{j+1/2,l}^{n+1/2})\\
       &-F(\rho_{j-1/2,r}^{n+1/2},R_{j+1/2,l}^{n+1/2})\bigg].
    \end{aligned}
\end{equation*}
The subindex $l$ and $r$ denote the reconstruction of $\rho$ and $R$ at $x_{\jph,l}$ and $x_{\jph,r}$, respectively:
\begin{eqnarray*}
\begin{aligned}\label{eq: mid values}
\rho_{j+1/2,l}^{n} &= \rho_j^{n} 
+  s_j^{n} \left(\frac{\Delta x}{2} + \Delta t c_{j+1/2}^-\right), \quad
\rho_{j+1/2,r}^{n} = \rho_{j+1}^{n} 
- s_{j+1}^{n} \left(\frac{\Delta x}{2} - \Delta t c_{j+1/2}^+\right).
\end{aligned}
\end{eqnarray*}
The midpoint values $\rho_{j+1/2,l/r}^{n+1/2}$ are approximated by using a Taylor expansion
\begin{align*}
   \rho_{j+1/2,l}^{n+1/2}&=\rho_{j+1/2,l}^n-\frac{\Delta t}{2}F_x(\rho_{j+1/2,l}^n,R_{j+1/2,l}^n),\\
   \rho_{j+1/2,r}^{n+1/2}&=\rho_{j+1/2,r}^n-\frac{\Delta t}{2}F_x(\rho_{j+1/2,r}^n,R_{j+1/2,r}^n).
\end{align*}
Here, the slopes $F_x(\rho_{j+1/2,l},R_{j+1/2,l})$ and $F_x(\rho_{j+1/2,r},R_{j+1/2,r})$ can be computed by the minmod limiter
\begin{small}
\begin{eqnarray*}
\begin{aligned}
&F_x(\rho_{j+1/2,l},R_{j+1/2,l})\\
&=\mm\left( \frac{F(\rho^{n}_{j+1/2,l},{R}^n_{j+1/2,l})-F(\rho_{{j-1/2,l}}^{n},{R}^n_{j-1/2,l})}{\Delta x-c_{j-1/2}^-\Delta t+ c_{j+1/2}^-\Delta t},\frac{F(\rho_{{j+3/2,l}}^{n},{R}^n_{j+3/2,l})-F(\rho_{{j+1/2,l}}^{n},{R}^n_{j+1/2,l})}{\Delta x-c_{j+1/2}^+\Delta t+ c_{j+3/2}^+\Delta t} \right),\\
&F_x(\rho_{j+1/2,r},R_{j+1/2,r})\\
&=\mm\left( \frac{F(\rho^{n}_{j+1/2,r},{R}^n_{j+1/2,r})-F(\rho_{{j-1/2,r}}^{n},{R}^n_{j-1/2,r})}{\Delta x-c_{j-1/2}^+\Delta t+ c_{j+1/2}^+\Delta t},\frac{F(\rho_{{j+3/2,r}}^{n},{R}^n_{j+3/2,r})-F(\rho_{{j+1/2,r}}^{n},{R}^n_{j+1/2,r})}{\Delta x-c_{j+1/2}^+\Delta t+ c_{j+3/2}^+\Delta t} \right).
\end{aligned}
\end{eqnarray*}
\end{small}So far the scheme is similar to the local fully-discrete scheme of \cite{kurganov2007reduction}. 
It remains to approximate the nonlocal terms.
For simplicity, we assume that the grid is chosen such that $N_\eta/\Delta x\in \N$.
Here, ${R}_{j+1/2,r/l}^{n+1/2}\approx (\wt\ast \rho)(t^{n+1/2},x^n_{j+1/2,r/l})$ is approximated by a Taylor series expansion, too
\begin{align*}
{R}(t^{n+1/2},x^n_{j+1/2,r/l})\approx
{R}(t^n,x^n_{j+1/2,r/l})+\frac{\Delta t}{2}{\partial_t R}(t^n,x^n_{j+1/2,r/l}).
\end{align*}
and the terms ${R}_{j+1/2,r/l},\,\partial_t R(t^n,x^n_{j+1/2,r/l})$ by midpoint quadrature rules: 
\begin{align*}%\label{conv : l}
R(t^n,x^n_{j+1/2,l})\approx&\int_{x_{j+1/2,l}}^{x_{j+1/2}}\omega_{\eta}(y-x_{j+1/2,l})\rho^{\Delta x}(t^n,y)dy \nonumber\\
+&\sum_{\rv{k}=0}^{N_\eta-2}\int_{x_{j+1/2+\rv{k}}}^{x_{j+3/2+\rv{k}}}\omega_{\eta}(y-x^n_{j+1/2,l})\rho^{\Delta x}(t^n,y)dy \nonumber\\
+&\int_{x_{j+N_\eta-1/2}}^{x_{j+1/2+N_\eta,l}}\omega_{\eta}(y-x^n_{j+1/2,l})\rho^{\Delta x}(t^n,y)dy,\nonumber\\
     =&c^-_{j+1/2}\Delta t\omega_{\eta}\left(-\frac{c^-_{j+1/2}\Delta t}{2}\right)\left(\rho_{j}^{n}+\frac{\Delta x+c^-_{j+1/2}\Delta t}{2}s^{n}_{j}\right) \nonumber\\
    +&\sum_{\rv{k}=0}^{N_\eta-2}\Delta x\omega_{\eta}((\rv{k}+1/2)\Delta x-c^-_{j+1/2}\Delta t)\rho_{j+\rv{k}+1}^{n} \\
    +&(\Delta x+c^-_{j+1/2}\Delta t)\omega_{\eta}\left(\frac{(2N_\eta-1)\Delta x-c^-_{j+1/2}\Delta t}{2}\right)\left(\rho_{j+N_\eta}^{n}+\frac{c^-_{j+1/2}\Delta t}{2} s_{j+N_\eta}^{n}\right).\nonumber
    \end{align*}
    Similarly, we obtain
\begin{align*}%\label{conv: r}
R(t^n,x^n_{j+1/2,r})  \approx&(\Delta x-c^+_{j+1/2}\Delta t)\omega_{\eta}\left(\frac{\Delta x-c^+_{j+1/2}\Delta t}{2}\right)\left(\rho_{{j+1}}^{n}+\frac{c^+_{j+1/2}\Delta t}{2}s^{n}_{{j+1}}\right) \nonumber\\
    +&\sum_{\rv{k}=0}^{N_\eta-2}\Delta x\omega_{\eta}((\rv{k}+3/2)\Delta x-c^+_{j+1/2}\Delta t)\rho_{j+\rv{k}+2}^{n}\\
    +&c^+_{j+1/2}\Delta t\omega_{\eta}\left(\frac{2N_\eta\Delta x-c^+_{j+1/2}\Delta t}{2}\right)\left(\rho_{j+N_\eta+1}^{n}-\frac{\Delta x-c^+_{j+1/2}\Delta t}{2}s_{j+N_\eta+1}^{n}\right). \nonumber
\end{align*}
To approximate $\partial_t R^k(t^n,x^n_{j+1/2,r/l})$ we use the conservation law \eqref{eq:conslaw} and a first-order quadrature rule. For the first and last cell we follow the same idea as discussed in \cite[Remark 2.9]{sudha2025convergence}, otherwise the midpoint rule. This keeps the computational cost low by using already computed values: 
\begin{small}
\begin{eqnarray*}
\begin{aligned}
\partial_t R(t^n,x^n_{j+1/2,r})=& -\int_{x^n_{j+1/2,r}}^{x^n_{j+1/2+N_\eta,r}}
\omega_{\eta}(y-x^n_{j+1/2,r})\partial_x F(\rho(t^n,y), (\wt\ast\rho)(t^n,y))dy,\\
    \approx &-(\Delta x-\Delta t c^+_{j+1/2})
    \omega_{\eta}\left(\frac{\Delta x-c^+_{j+1/2}\Delta t}{2}\right)F_x(\rho^{n}_{j+1/2,r},{R}^{n}_{j+1/2,r})\\
    &-\sum_{\rv{k}=0}^{N_\eta-2}\Delta x\omega_{\eta}((\rv{k}+3/2)\Delta x-c^+_{j+1/2}\Delta t )F_x(\rho^{n}_{j+3/2+\rv{k},r},{R}^{n}_{j+3/2+\rv{k},r})\\
    &-c^+_{j+1/2}\Delta t\omega_{\eta}\left(\frac{2N_\eta\Delta x-c^+_{j+1/2}\Delta t}{2}\right)F_x(\rho^{n}_{j+1/2+N_\eta,r},{R}^{n}_{j+1/2+N_\eta,r}),\\
\partial_t R(t^n,x^n_{j+1/2,l})
\approx& c^-_{j+1/2}\Delta t\omega_{\eta}\left(\frac{-c^-_{j+1/2}\Delta t}{2}\right) F_x(\rho^{n}_{j+1/2,l},R^{n}_{j+1/2,l})\\
    &-\sum_{\rv{k}=0}^{N_\eta-2}
    \Delta x\omega_{\eta}((\rv{k}+1/2)\Delta x-c^-_{j+1/2}\Delta t) F_x(\rho^{n}_{j+3/2+\rv{k},l},{R}^{n}_{j+3/2+\rv{k},l})\\
    &-(\Delta x+c^-_{j+1/2}\Delta t)\omega_{\eta}\left(\frac{(2N_\eta-1)\Delta x-c^-_{j+1/2}\Delta t}{2}\right)
    F_x(\rho^{n}_{j+1/2+N_\eta,l},{R}^{n}_{j+1/2+N_\eta,l}).
\end{aligned}
\end{eqnarray*}    
\end{small}
We note that the spatial derivatives of the flux are already computed via the minmod limiter.
Finally, it remains to project the approximated solution back on the original grid, which can be achieved by using another piecewise linear reconstruction 
\begin{equation*}
    \begin{aligned}
      {w}^{\Delta x}(t^{n+1},x)=&\sum_{j\in\Z} \left( w_{j+1/2}^{n+1}+s_\jph^{n+1}\left(x-\frac{x_{\jph,l}+x_{\jph,r}}{2}\right)\right)\chi_{[x_{j+1/2,r}^n,x_{j+1/2,l}^n]}(x)\\ &+w_j^{n+1}\chi_{[x_{j-1/2,r}^n,x_{j+1/2,l}^n]} (x).
    \end{aligned}
\end{equation*}
Here, the slopes are computed similar to \cite{kurganov2007reduction}
\begin{align*}
    s_\jph^{n+1}=\mm\left(\frac{w_\jph^{n+1}-\rho_{\jph,l}^{n+1}}{(c_\jph^+-c_\jph^-)\Delta t/2},\frac{\rho_{\jph,r}^{n+1}-w_\jph^{n+1}}{(c_\jph^+-c_\jph^-)\Delta t/2}\right),
\end{align*}
where $\rho_{\jph,l/r}^{n+1}$ are obtained by another Taylor expansion. In addition, monotonicity is enforced if the values are not between $w_j^{n+1}$ and $w_\jph^{n+1}$. For more details we refer to \cite{kurganov2007reduction}.
Integrating the piecewise linear function over the original grid results in the scheme
\begin{equation}\label{eq 3}
    \begin{aligned}
        \rho_j^{n+1}
        =&w_{j}^{n+1}+\frac{\Delta t}{\Delta x}\bigg[c_{j-1/2}^{+}(w_{j-1/2}^{n+1}-w_j^{n+1})-c_{j+1/2}^{-}(w_{j+1/2}^{n+1}-w_j^{n+1})\bigg]\\
        &+\frac{\Delta t^2}{2 \Delta x}\bigg[s_\jph^{n+1}c_\jph^+c_\jph^--s_\jmh^{n+1}c_\jmh^+c_\jmh^-\bigg].
    \end{aligned}
\end{equation}
The disadvantage of the fully-discrete scheme for nonlocal conservation laws is that it has rather high computational cost: in total we need to compute five nonlocal terms: one for the speed estimates, four for the evolution of the scheme. Additionally, the weights for the latter ones need to be recomputed every time step.
All this causes a computational bottleneck.
In addition, the fully-discrete form does not admit a straightforward extension to higher orders than two (extensions to suitable systems and balance laws are possible). Hence, usually a semi-discrete form of the  fully-discrete scheme is derived by considering the limit of the time discretization.

\subsection{The central-upwind flux}
The CU flux the nonlocal problem \eqref{eq:conslaw}--\eqref{eq:Fsimple} can be derived following the same procedure as in \cite{kurganov2007reduction}. 
In particular, we consider $\lim_{\Delta t \rightarrow 0}\frac{\rho_j^{n+1}-\rho_j^n}{\Delta t^n}$ which results in the semi-discrete CU scheme:
\begin{align}\label{eq:semidiscrete}
    \frac{d}{dt}\bar{\rho}_j(t)=-\frac{\mathcal{F}_{j+1/2}(t)-\mathcal{F}_{j-1/2}(t)}{\Delta x},
\end{align}
where $\bar \rho_j(t)$ denotes the cell-average at time $t$ and $\mathcal{F}_{j+1/2}$ is given by
\begin{align}\label{eq:upwindflux}
    \mathcal{F}_{j+1/2}(t)=&\frac{c_{j+1/2}^{+}F(\rho_{j+1/2}^{-},R_{j+1/2})-c_{j+1/2}^{-}F(\rho_{j+1/2}^{+},R_{j+1/2})}{c_{j+1/2}^{+}-c_{j+1/2}^{-}}\\
    &\nonumber+\frac{c_{j+1/2}^{+}c_{j+1/2}^{-}}{c_{j+1/2}^{+}-c_{j+1/2}^{-}}\bigg[\rho_{j+1/2}^{+}-\rho_{j+1/2}^{-}-d_{j+1/2}\bigg].
\end{align}
Here, $d_{j+1/2}$ is the "built-in" anti-diffusion term given by
\begin{align*}d_{j+1/2}=&\mm\bigg(\rho_{j+1/2}^{+}-\rho_{j+1/2}^{*},\,\rho_{j+1/2}^{*}-\rho_{j+1/2}^{-}\bigg),
\end{align*}
with $\rho_{j+1/2}^{*}$ obtained by
\begin{align*}
   \rho_{j+1/2}^{*}=\frac{c_{j+1/2}^{+}\rho_{j+1/2}^{+}-c_{j+1/2}^{-}\rho_{j+1/2}^{-}-\left(F(\rho_{j+1/2}^{+},R_{j+1/2})-F(\rho_{j+1/2}^{-},R_{j+1/2})\right)}{c_{j+1/2}^{+}-c_{j+1/2}^{-}}.
\end{align*}
The flux is very similar to the one derived in \cite{kurganov2009non}: The main differences are that we consider a more general flux function and consider the nonlocal term as fixed at the cell interface to get an estimate on the speed while in \cite{kurganov2009non} an upper bound on the nonlocal term is used to estimate the speed which impacts particularly the anti-diffusion term.\\
As aforementioned, the main advantage of the CU schemes is that it can be used as a black-box solver for general nonlocal conservation laws, while keeping the numerical diffusion much lower as with a Lax-Friedrichs-type scheme.
Due to the semi-discrete formulation, we can use the CU flux to obtain higher-order approximations by using appropriate spatial reconstructions and Runge-Kutta methods in time.

\subsection{Convergence of the central-upwind scheme}
We start by discussing the first-order CU scheme. 
Therefore, the semi-discrete form is solved by an explicit Euler scheme in time.
The spatial reconstruction is given by a piecewise constant function $\dr$ such that $\dr(t,x)=\rho_j^n$ for $(t,x)\in[t^n,t^{n+1})\times [x_\jmh,x_\jph)$ with $\rho_j^n$ approximating the cell-averages, see \eqref{eq:cellavg}.
Hence, we have
\begin{equation*}
\begin{aligned}
&\rho_{j+1/2}^+={\rho}_{j+1}^n,\quad\rho_{j+1/2}^-={\rho}_j^n.
\end{aligned}
\end{equation*}
We approximate the initial data by
\begin{align*}
    \rho_j^0=\frac{1}{\Dx} \int_{x_\jmh}^{x_\jph} \rho_0(x) dx,\quad j\in\Z.
\end{align*}
Over a piecewise constant function the approximation of the nonlocal term in \eqref{eq:approximatedNLterm} is exact (with $s_{j+\Ne+1}^n=0$\rv{)}.
Since we are considering a first-order approximation, we can set $\gamma_{\Ne+1}=0$, too, which results in an even simpler approximation.

First, we note that if $g'(\rho)\geq 0$ for $\rho\in I$, the CU flux reduces to the widely used Godunov/Upwind-type approximations for nonlocal conservation laws, which where originally introduced in \cite{friedrich2018godunov}, i.e. the CU flux is simply given by
\[\mathcal{F}_{j+1/2}(t)=g(\rho_j) v(R_{j+1/2}).\]

Let us denote by $\|f\|:=\|f\|_{L^\infty(I,\R)}$.
Then, we have the following result, for the general case ($g$ nonlinear, but $g$ convex/concave),
\begin{theorem}\label{thm}
    Using an explicit Euler scheme for the semi-discrete scheme \eqref{eq:semidiscrete} with the central-upwind flux as in \eqref{eq:upwindflux}, the resulting numerical scheme converges, under the CFL condition 
    \begin{equation}\label{eq:CFLCU}
\frac{\Delta t}{\Delta x}\leq \frac{1}{(2\norm{g'}\rv{\rho_{\max}}+\norm{g}) \norm{v'} \gamma_0+4\norm{g'}\norm{v}}
\end{equation}
    for $(\Dt,\Dx)\to 0$ to the unique weak entropy solution of \eqref{eq:conslaw}--\eqref{eq:Fsimple} in the sense of \cite[Definition 2.2]{friedrich2023numerical} and satisfies the maximum principle $\rho_j^n\in I$ for $j\in \Z$ and $n\in\N$.
\end{theorem}
\rv{
\begin{remark}
    We note that the CFL condition \eqref{eq:CFLCU} is more restrictive than the one given by \eqref{eq:CFL0}, since for the considered flux function it holds that 
    \[\max_j \left( \max \left( c^+_{j+\frac{1}{2}}, -c^-_{j+\frac{1}{2}} \right) \right)\leq \norm{g'}\norm{v}.\]
\end{remark}
}
\begin{proof}
We set $\mathcal{F}_\jph^n:=\mathcal{F}_\jph(t^n)$. In addition, it depends on the inputs $(a)$ $\rho_\jph^-=\rho_j^n$ and $(b)$ $\rho_\jph^+=\rho_\jpo^n$, such that we write $\mathcal{F}_\jph^n(a,b)$.
A numerical flux function  $\mathcal{F}^n_\jph$ converges to the weak entropy solution of \eqref{eq:conslaw}--\eqref{eq:Fsimple} and satisfies the maximum principle, if the following conditions are fulfilled, cf. \cite[Definition 3.1--3.2,Theorem 4.1]{friedrich2023numerical}:
\begin{enumerate}
\item[\emph{(i)}] $\mathcal{F}^n_\jph(a,b)=G(a,b)v(R_{j+1/2}^n)$
with $G$ only depending on $g,\ \rho_{\min},\ \rho_{\max}$. 
\item[\emph{(ii)}] Consistency: $\mathcal{F}^n_\jph(\rho,\rho)=g(\rho)v(R_{j+1/2}^n)\ \forall \rho\in [\rho_{\min},\rho_{\max}]$.
\item[\emph{(iii)}] The map $(a,b)\mapsto \mathcal{F}^n_\jph(a,b)=G(a,b)v(R_{j+1/2}^n)$ from $[\rho_{\min},\rho_{\max}]^2$ to $\R$ is non-decreasing with respect to $a$ and non-increasing with respect to $b$.
\item[\emph{(iv)}] $\mathcal{F}^n_\jph(a,b)$ satisfies 
\begin{align*}
    \abs{\mathcal{F}^n_\jph(a,b)-\mathcal{F}^n_\jph(b,b)}\leq L_1 \abs{a-b}\quad\text{and}\quad\abs{\mathcal{F}^n_\jph(a,b)-\mathcal{F}^n_\jph(a,a)}\leq L_2\abs{a-b}
\end{align*}
    for $(a,b)\in [\rho_{\min},\rho_{\max}]^2$ and every $R_{j+1/2}^n\in I$.
 \item[\emph{(v)}] The CFL condition 
\begin{equation*}
\frac{\Dt}{\Dx}\leq \frac{1}{\norm{G} \norm{v'} \gamma_0+L_1+L_2}
\end{equation*}
holds.
\end{enumerate}
Hence, we need to show that the flux function \eqref{eq:upwindflux} fulfills the conditions \emph{(i)--(iv)} and derive the CFL condition.
Let us consider the speeds $c_{j+1/2}^\pm$.
Due to the flux \eqref{eq:Fsimple}, $g$ being convex/concave and $v: I\to\R_+$, we can write
\[c_{j+1/2}^+(a,b)=\tilde c_{j+1/2}^+(a,b)v(R_{j+1/2}^n)\quad \text{with}\quad \tilde c_{j+1/2}^+(a,b):=\max\{g'(a),g'(b),0\}.\]
We can define $\tilde c_{j+1/2}^-(a,b)$ in the same way.
This simplifies the nonlocal numerical flux function such that we can rewrite it as
\[\mathcal{F}^n_\jph(a,b)=G(a,b)v(R_{j+1/2}^n),\]
where $G$ is the usual (local) CU flux applied on $g$.
This proves \emph{(i)}.\\
The condition \emph{(ii)} follows immediately and \emph{(iii)} follows from the monotonicity of the CU scheme, see \cite[Remark 2.2]{cui2024bound} and \cite{bryson2005semi}.\\
It remains to prove $(iv)$. Let us first observe that by \emph{(i)} and \emph{(ii)} the problem reduces to proving $|G(a,b)-g(b)|\leq \tilde L_1|a-b|$ and $|G(a,b)-g(a)|\leq \tilde L_2|a-b|$. In the following we omit for readability the dependence of $c_{j+1/2}^\pm$ on $a,b$ and estimate
\begin{align*}
    b-\rho^*_{j+1/2}&=\frac{(\tilde c_{j+1/2}^+-\tilde c_{j+1/2}^-)b-\tilde c_{j+1/2}^+ b+\tilde c_{j+1/2}^-a +g(b)-g(a)}{\tilde c_{j+1/2}^+-\tilde c_{j+1/2}^-}\\
    &\leq \frac{-\tilde c_{j+1/2}^-+\norm{g'}}{c_{j+1/2}^+-\tilde c_{j+1/2}^-} |a-b|.
\end{align*}
We obtain a similar estimate for $\rho^*_{j+1/2}-a$, such that we can estimate
$$d_{j+1/2}\leq \frac{\max\{c_{j+1/2}^+,\ -\tilde c_{j+1/2}^-\}+\norm{g'}}{c_{j+1/2}^+-\tilde c_{j+1/2}^-}|a-b|. $$
Now, we look at the remaining part of $|G(a,b)-g(b)|$, i.e.
\begin{align*}
    &\frac{\tilde c_{j+1/2}^+g(a)-\tilde c_{j+1/2}^-g(b)+\tilde c_{j+1/2}^+\tilde c_{j+1/2}^-(b-a)-(\tilde c_{j+1/2}^+-\tilde c_{j+1/2}^-)g(b)}{\tilde c_{j+1/2}^+-\tilde c_{j+1/2}^-}\\
    &=   \frac{\tilde c_{j+1/2}^+(g(a)-g(b))+\tilde c_{j+1/2}^+\tilde c_{j+1/2}^-(b-a)}{\tilde c_{j+1/2}^+-\tilde c_{j+1/2}^-}\leq \frac{\tilde c_{j+1/2}^+\norm{g'}-\tilde c_{j+1/2}^+\tilde c_{j+1/2}^-}{\tilde c_{j+1/2}^+-\tilde c_{j+1/2}^-}|a-b|
\end{align*}
Combining both estimates and using $|\tilde c_{j+1/2}^\pm|\leq \norm{g'}$ gives 
\begin{small}
\begin{align*}
    &|G(a,b)-g(b)|\\
    &\leq \frac{(\tilde c_{j+1/2}^+\norm{g'}-\tilde c_{j+1/2}^+\tilde c_{j+1/2}^-)(\tilde c_{j+1/2}^+-\tilde c_{j+1/2}^-)-\tilde c_{j+1/2}^+\tilde c_{j+1/2}^-(\max\{c_{j+1/2}^+,\ -\tilde c_{j+1/2}^-\}+\norm{g'})}{(\tilde c_{j+1/2}^+-\tilde c_{j+1/2}^-)^2}|a-b|\\
    &\leq \frac{(\tilde c_{j+1/2}^+)^2-\tilde c_{j+1/2}^+\tilde c_{j+1/2}^- +(\tilde c_{j+1/2}^+)^2+ (\tilde c_{j+1/2}^-)^2-2\tilde c_{j+1/2}^+ \tilde c_{j+1/2}^-}{(\tilde c_{j+1/2}^+-\tilde c_{j+1/2}^-)^2}\norm{g'}|a-b|\\
    &\leq \frac{(\tilde c_{j+1/2}^+)^2-2\tilde c_{j+1/2}^+\tilde c_{j+1/2}^- + (\tilde c_{j+1/2}^-)^2}{(\tilde c_{j+1/2}^+-\tilde c_{j+1/2}^-)^2}2\norm{g'}|a-b|=2\norm{g'}|a-b|.
\end{align*}
\end{small}
 The same estimate can be obtained for the second inequality analogously.\\
 For the CFL condition we can use the estimates above and $|G(a,b)|\leq |G(a,b)-g(b)|+g(b)$ to obtain \eqref{eq:CFLCU}.
\end{proof}
Now, we consider a second-order accurate CU scheme. 
Therefore, we use the piecewise linear reconstruction \eqref{eq:piecewise} in space and use the two-stage second-order accurate strong-stability-preserving (SSP) Runge-Kutta scheme in time, see \cite{GottliebShu1998,gd2023convergence}.
We note that this second-order scheme only requires the computation of two nonlocal terms and hence the computational cost are significant\rv{ly} lower in comparison to the KT scheme.
\rv{We get the following result}
\begin{corollary}
    \rv{Let $g’\geq 0$.} Using the second-order SSP Runge-Kutta scheme for the semi-discrete scheme \eqref{eq:semidiscrete} with the central-upwind flux as in \eqref{eq:upwindflux} and the piecewise linear spatial reconstruction \eqref{eq:piecewise} with the modification of the slope limiters as in \cite[eq. (5.3)]{gd2023convergence}, the resulting numerical scheme converges, under the CFL condition 
    \begin{equation}\label{eq:CFLCU2}
\frac{\Delta t}{\Delta x}\leq \frac{1}{2(\norm{g} \norm{v'} \gamma_0+\norm{g'}\norm{v})}
\end{equation}
    for $(\Dt,\Dx)\to 0$ to the unique weak entropy solution of \eqref{eq:conslaw}--\eqref{eq:Fsimple} in the sense of \cite[Definition 2.2]{friedrich2023numerical} and satisfies the maximum principle $\rho_j^n\in I$ for $j\in \Z$ and $n\in\N$.
\end{corollary}
\begin{proof}
    For $g’\geq 0$ the CU flux coincides with a Godunov-type flux. The convergence of a second-order Godunov-type flux with the nonlocal term approximated by the trapezoidal rule is proven in \cite{gd2023convergence}.
    Here, we use a different approximation, i.e. \eqref{eq:approximatedNLterm}.
    However, for a grid chosen such that $\eta/\Delta x\in \N$, all the necessary estimates on the differences of the nonlocal terms are established in \cite{friedrich2023numerical} by using additionally $\rho_j^n=(\rho_\jph^-+\rho_\jmh^+)/2$.
    These estimates can also be extended to the general case $\eta/\Delta x\notin \N$.
    Hence, the convergence result follows from the proofs in \cite{gd2023convergence}.
\end{proof}
\begin{remark}
    For the general case of $g$ being convex/concave the convergence of the second-order scheme can also be proven. 
    However, it is not as straightforward as for the case $g'\geq 0$ and, therefore, it is beyond the scope of this work.
    While the proof of a discrete maximum principle follows immediately from the proof of Theorem \ref{thm} and \cite[Theorem 3.1]{gd2023convergence},
    the derivation of a BV bound as in \cite[Proposition 4.5]{gd2023convergence} needs to be studied in detailed. We expect that following the proofs of \cite[Proposition 4.5]{gd2023convergence} and \cite[Lemma 3.2]{friedrich2023numerical} such a bound can be obtained which allows to follow then the same steps as in \cite{gd2023convergence} to prove the convergence.
\end{remark}
\begin{remark}
    High-order approaches (order three and above) using the CU flux in the case $g'\geq 0$ have been (although not intentionally) studied in \cite{friedrich2019maximum}.
\end{remark}
\subsection{Extension to \rv{weakly coupled systems} of nonlocal balance laws}
We close this section by commenting on the extension of the CU and KT scheme to systems of nonlocal conservation and balance laws.
For a general nonlocal system, i.e. 
\[\partial_t \brho +\partial_x F(\brho,\omega \ast \brho)=0,\]
where $\brho(t,x)\in\mathbb{R}^{N},\ (\omega\ast \brho) (t,x) \in\mathbb{R}^{m},\,F:\R^N\times \R^m\to \mathbb{R}^{N}$, it is not straightforward to derive numerical schemes (even in the local case).
However, central schemes can be used.
In particular, the derivation of the KT scheme and hence, the CU scheme remains valid for systems of nonlocal conservation laws.
The challenging part is the estimation of the speeds.
As in the local case, the speeds need to be estimated by using the Jacobian of $F(\brho,\omega \ast \brho)$.
One way to estimate them could be to fix the nonlocal term at the cell-interface and compute the Jacobian of $F(\rho,R(t,x_{j+1/2}))$ concerning $\rho$.\\
However to best of the authors' knowledge, in current research systems of nonlocal conservation laws are typically only coupled through the nonlocal term (or possible source terms).
Hence, the Jacobian (for a fixed nonlocal term) becomes a diagonal matrix and the speeds can be estimated. 
In particular, due to the weak coupling we can use a component-wise approach.
To illustrate this more clearly, we consider the following system similar to \cite{aggarwal2015nonlocal,aggarwal2024well}
\begin{equation}\label{eq:system}
    \partial_t \rho_k +\partial_x F_k(\rho_k,\omega \ast \brho)=0,
\end{equation}
for $k=1,\dots, N$. For the precise assumptions on the parameters we refer to \cite{aggarwal2015nonlocal,aggarwal2024well}, but let us additionally assume that each $F_k$ is either convex or concave in its first component.
The simple coupling allows us to apply the first-order CU scheme component-wise and to estimate the speeds for each $k=1,\dots, N$ separately:
\begin{equation*}
\begin{aligned}
c_{k,j+\frac12}^{+} &=
\max\left(
\frac{\partial F_k(\rho_{k,j+\frac12}^-,(\omega \ast \brho)(t,x_\jph))}{\partial \rho},
\frac{\partial F_k(\rho_{k,j+\frac12}^+,(\omega \ast \brho)(t,x_\jph))}{\partial \rho},
0
\right), \\
c_{k,j+\frac12}^{-} &=
\min\left(
\frac{\partial F_k(\rho_{k,j+\frac12}^-,(\omega \ast \brho)(t,x_\jph))}{\partial \rho},
\frac{\partial F_k(\rho_{k,j+\frac12}^+,(\omega \ast \brho)(t,x_\jph))}{\partial \rho},
0
\right).
\end{aligned}
\end{equation*}
Similar to the scalar case if each $F_k$ is a linear function in its first argument, the component-wise first-order CU scheme coincides with the Godunov-type scheme (simplifying to a Upwind-type scheme in the linear case).
The resulting scheme was analyzed in \cite{aggarwal2024well} and the convergence to the unique weak (entropy) solution is proven.\\
\rv{There are several systems of nonlocal conservation laws in applications, which are linear in the local variable.
Most of these systems model traffic flow. Examples include multi-class traffic flow \cite{chiarello2019multiclass}, bidirectional traffic flow \cite{chiarello2021non} as well as the nonlocal Keyfitz-Kranzer model \cite{aggarwal2015nonlocal,aggarwal2024well}.
Here, the authors employed an Upwind-type scheme and implicitly the CU-type scheme to demonstrate the convergence of the scheme and the existence of solutions.}
Another system with a linear dependence in the local convective part are multi-commodity supply chain models \cite{gugat2016analysis,keimer2018analysis}.
We note that a proof of convergence for the systems \eqref{eq:system} of a first-order CU scheme remains open, but we expect that it can be obtained similar to \cite{aggarwal2024well}.
Note that similar to the CU scheme the KT scheme can be directly applied to \eqref{eq:system}.
\\
Now let us consider systems of nonlocal balance laws, similar to \cite{sudha2025convergence}
\begin{equation*}
    \partial_t \rho_k +\partial_x F_k(\rho_k,\omega \ast \brho)=S_k(\brho,\omega \ast \brho),
\end{equation*}
for $k=1,\dots,N$. We refer to \cite{sudha2025convergence} for details on the parameters.
Different approaches are possible to deal with the source terms.
One way is to treat it by a half-time step. Here, we refer to \cite{friedrich2020nonlocal} which proves the convergence of the CU scheme (coinciding once more with the Godunov-type scheme) for a multilane traffic model.
Alternatively, the source can be added immediately without a half-step.
For the second-order accurate CU scheme, we can use the trapezoidal rule to approximate the source term and for the KT scheme we need to follow similar approaches as in \cite{sudha2025convergence}. In particular, whenever the balance law is used after the Taylor expansions, the source terms need to be included. For a second-order accurate approximation, the trapezoidal rule should be used in space to avoid additional computational effort. For more details we refer to \cite{sudha2025convergence}. Here, a similar scheme (without dividing the cells into smooth and non-smooth parts) is discussed.\\
In addition, we can apply the CU type scheme to systems of nonlocal balance laws, which can have a mixed local and nonlocal flux and e.g. an additional coupling by the source term. Examples of such systems are the nonlocal Euler equations in \cite{bhatnagar2021well} or the two-lane bidirectional traffic model with overtaking from \cite{contreras2025two}.

\section{Numerical examples}\label{sec:3}
We consider two examples: the scalar case which we discussed in detail above and the case of a system of balance laws. All reference solutions in this section are computed by a second-order Nessyahu-Tadmor scheme and $\Delta x= \frac{1}{20}\cdot 2^{-9}$. The convergence of this scheme for systems of nonlocal balance laws is proven in \cite{sudha2025convergence}. In all examples the kernel has compact support on $[0,\eta]$ with $\eta=0.2$ and is given by the quadratic function $\omega_{\eta}(x)=3\frac{\eta^{2}-x^2}{2\eta^3}$. In addition, we consider the interval $[-1,1]$ with periodic boundary conditions.
\rv{The time step sizes are chosen according to the minimum of the different CFL conditions in the following. 
We use the upper bounds of the CFL conditions \eqref{eq:CFLCU} and \eqref{eq:CFLCU2}. For KT scheme, the CFL condition derived for the second-order Nessyahu-Tadmor scheme in \cite[eq. (2.14)]{sudha2025convergence} provides stable and maximum-principle-satisfying solutions in all tests: 
\[\Delta t\leq \frac{\sqrt{2}-1}{2}\frac{\Delta x}{\norm{g'}\norm{v}}.\]
This is the most restrictive CFL condition, particularly when the step sizes are small.}
\subsection{Scalar conservation law}
The Arrhenius traffic flow model is given by $g(\rho)=\rho(1-\rho)$ and $v(\rho)=\exp(-\rho)$. Apart from the CU scheme (first- and second-order) and the \rv{KT} scheme, we consider a Godunov-type scheme following \cite{friedrich2023numerical}.
First, we test the numerical convergence rate in the $L^1$-norm at $T=0.15$ for smooth initial data
\begin{align}
    \rho_0(x)=0.5+0.4\sin(\pi x),
\end{align}
and $\Dx=\frac{1}{20}\cdot2^{-n}$ for $n=0,\dots,5$. The computed convergence rates in Table \ref{tab:Traffic cu and GD} show that all schemes achieve the expected order of convergence. We note that there are no significant differences between the CU (first-order) and Godunov type flux and similarly between the second-order CU and KT scheme. 

\begin{figure}
     \centering
     \setlength{\fwidth}{0.8\linewidth}
     \input{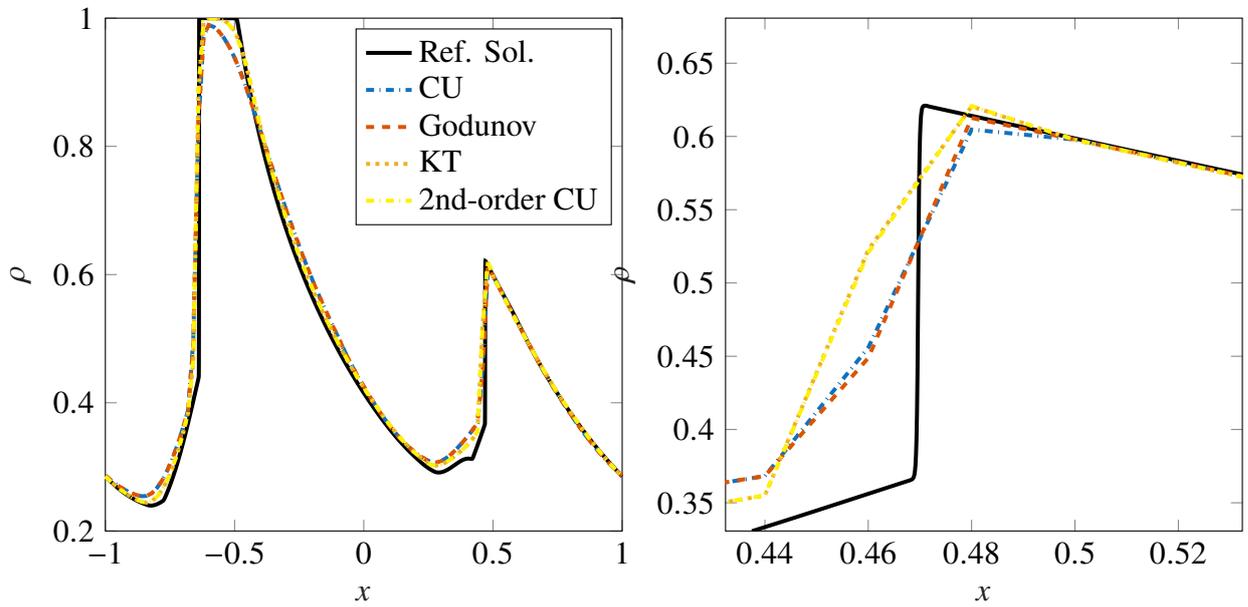}
     \caption{Approximate solutions of the CU (first- and second-order), Godunov and KT scheme for a discontinuous initial data at $T=1$ with $\Delta x=0.02$. The right panel shows a zoom into the solution.}\label{fig:Arrhenius}
 \end{figure}

 \begin{table}[hbt!]
 \centering
 \caption{Comparison of L$^1$ errors and convergence orders for CU (first- and second-order), Godunov and KT schemes.}
\label{tab:Traffic cu and GD}
     \begin{tabular}{|c | c  c | c  c | c c|c c|}
     \hline
         &\multicolumn{2}{c|}{CU} & \multicolumn{2}{c|}{Godunov} & \multicolumn{2}{c|}{KT}& \multicolumn{2}{c|}{2nd-order CU} \\
$n$&$L^1$-error&c.r. &$L^1$-error&c.r. &$L^1$-error&c.r.&$L^1$-error&c.r. \\
\hline
0&7.83e-03&-&7.86e-03&-&1.45e-03&-&1.54e-03&-\\
1&4.12e-03&0.93&4.14e-03&0.92&3.97e-04&1.87&4.33e-04&1.83\\
2&2.08e-03&0.99&2.08e-03&1.00&1.06e-04&1.91&1.15e-04&1.92\\
3&1.04e-03&0.99&1.04e-03&0.99&2.80e-05&1.91&2.98e-05&1.94\\
4&5.23e-04&1.00&5.23e-04&1.00&7.48e-06&1.90&7.80e-06&1.93\\
5&2.62e-04&1.00&2.62e-04&1.00&1.92e-06&1.96&1.92e-06&2.02\\
 \hline
    \end{tabular} 
\end{table}
However, for a discontinuous initial data as, e.g.
\[\rho_0(x)=\begin{cases}
    1,&x \in[-0.5,0],\\
    0.8,&x \in[0.5,0.75],\\
    0.2,&\text{else},
\end{cases}\]
more significant differences between the first-order schemes can be obtained. The approximate solutions at $T=1$ and $\Delta x=0.02$ are displayed in Figure \ref{fig:Arrhenius}. In particular, the zoom on the right shows that the Godunov type scheme provides a slightly more accurate approximation than the first-order CU scheme. Again differences between the KT and second-order CU scheme are not visible.

\subsection{Systems of balance laws}
Now, we consider the nonlocal multilane model \cite{bayen2022multilane,friedrich2020nonlocal} given by
\begin{align*}
    \partial_t \rho^1+ \partial_x (\rho^1 v(\omega_\eta\ast \rho^1))&=-S(\rho^1,\rho^2,\omega_\eta\ast \rho^1,\omega_\eta\ast \rho^2),\\
    \partial_t \rho^2+ \partial_x (\rho^2 v(\omega_\eta\ast \rho^2))&=S(\rho^1,\rho^2,\omega_\eta\ast \rho^1,\omega_\eta\ast \rho^2).
\end{align*}
The source term is given similar to \cite{friedrich2020nonlocal} as
\[S(\rho^1,\rho^2,R^1,R^2)=(v(R^2)-v(R^1))\begin{cases}
    \rho^1(1-\rho^2),&\text{if }v(R^2)\geq v(R^1),\\
    \rho^2(1-\rho^1),&\text{if }v(R^2)< v(R^1).
\end{cases}\]
As aforementioned, the CU scheme coincides with the Godunov-type scheme and the convergence is proven in \cite{friedrich2020nonlocal}. The speed function is chosen as $v(\rho)=1-\rho^2$ and we start by testing the convergence for the smooth initial data
\[\rho_0^1(x)=0.5+0.5\sin(\pi x),\quad \rho_0^2(x)=0.25+0.25\cos(2\pi x),\]
at the final time \rv{$T=0.5$}. All schemes attain the expected order of convergence and the differences between the second-order CU scheme and the KT scheme are small, see Table \ref{tab:multilane}.
\rv{Additionally, we display upper bounds on computational time. These demonstrate that the second-order CU scheme is more efficient than the KT scheme. 
Note that we used the same time step sizes for both schemes.
In particular, larger time steps are possible with the second-order CU scheme, making it even more efficient.
Results obtained by selecting the time step according to the upper bound provided by  \eqref{eq:CFLCU2} are presented in Table \ref{tab:multilaneCU}.
The computational costs decrease further while the errors increase only slightly.}
 \begin{table}[hbt!]
 \centering
 \caption{Comparison of L$^1$ errors\rv{, convergence orders and computational time} for CU, a second-order CU and KT schemes. 
 \rv{ The computational time is given as an upper bound in seconds.}}
\label{tab:multilane}
     \begin{tabular}{|c | c c c | c c c | c c c|}
     \hline
         &\multicolumn{3}{c|}{CU} & \multicolumn{3}{c|}{2nd order CU} & \multicolumn{3}{c|}{KT} \\
$n$&$L^1$-error&c.r.&\rv{c.t. (s)}&$L^1$-error&c.r.&\rv{c.t. (s)} &$L^1$-error&c.r.&\rv{c.t. (s)} \\
\hline
 0&1.30e-01&-& \rv{$<0.01$}&4.24e-02&-& \rv{$<0.01$}&4.28e-02&-& \rv{$<0.01$}\\
1&7.55e-02&0.79& \rv{$<0.01$}&1.61e-02&1.39& \rv{$<0.01$}&1.61e-02&1.41& \rv{$<0.02$}\\
2&4.15e-02&0.87& \rv{$<0.01$}&4.97e-03&1.70& \rv{$<0.02$}&4.92e-03&1.71& \rv{$<0.05$}\\
3&2.21e-02&0.91& \rv{$<0.02$}&1.39e-03&1.84& \rv{$<0.05$}&1.36e-03&1.85& \rv{$<0.20$}\\
4&1.15e-02&0.94& \rv{$<0.05$}&3.69e-04&1.91& \rv{$<0.10$}&3.63e-04&1.91& \rv{$<0.80$}\\
5&5.87e-03&0.97& \rv{$<0.20$}&9.42e-05&1.97& \rv{$<0.35$}&9.24e-05&1.97& \rv{$<4.65$}\\
 \hline
    \end{tabular} 
\end{table}

\begin{table}
    [hbt!]
 \centering
 \caption{\rv{Comparison of L$^1$ errors, convergence orders and computational time for a second-order CU scheme with a less restrictive time step. The computational time is given as an upper bound in seconds.}}
\label{tab:multilaneCU}
     \begin{tabular}{|c | c c c |}
     \hline
          & \multicolumn{3}{c|}{2nd order CU} \\
$n$&$L^1$-error&c.r.&{c.t. (s)} \\
\hline
0&4.24e-02&-&$<0.01$\\
1&1.64e-02&1.37&$<0.01$\\
2&5.09e-03&1.69&$<0.01$\\
3&1.42e-03&1.84&$<0.02$\\
4&3.80e-04&1.90&$<0.05$\\
5&9.75e-05&1.96&$<0.13$\\
 \hline
    \end{tabular} 
\end{table}

\begin{figure}
     \centering
     \setlength{\fwidth}{0.8\linewidth}
     \input{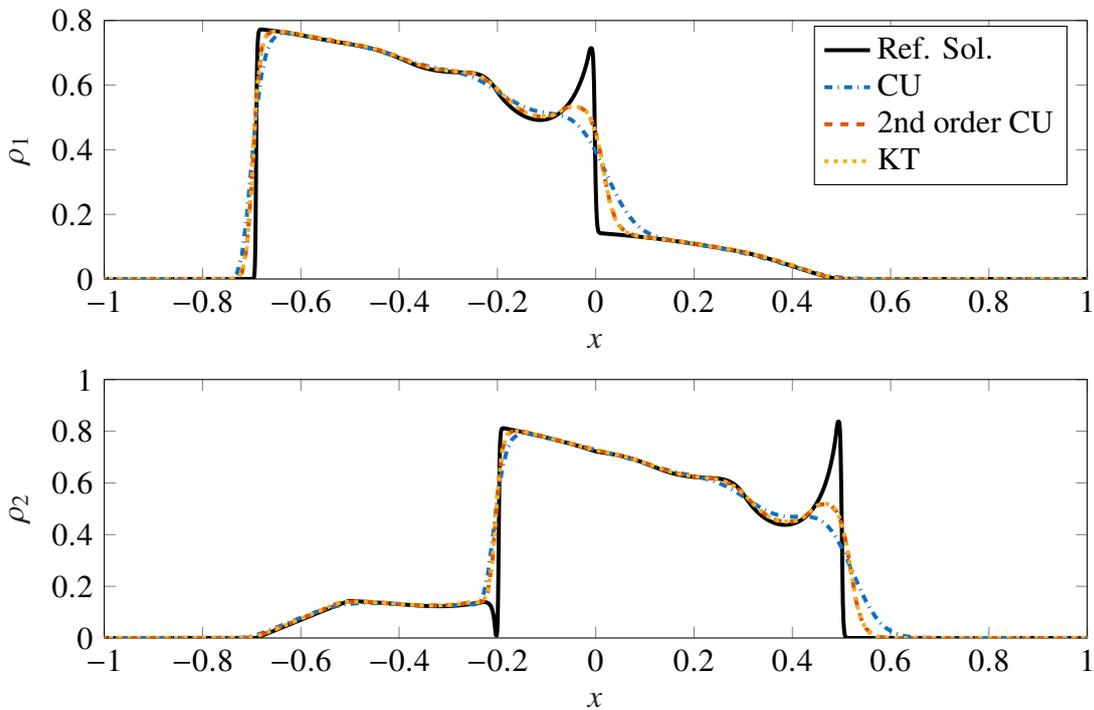}
     \caption{Approximate solutions of the CU, second-order CU and KT scheme for a discontinuous initial data at $T=0.25$ with $\Delta x=0.01$.}\label{fig:Multilane}
 \end{figure}

The same effect can be seen for the approximate solutions for the discontinuous initial data 
\[\rho_{0,1}(x)=1\chi_{[0,0.5]}(x),\quad \rho_{0,2}(x)=1\chi_{[0.5,1]}(x)\]
in Figure \ref{fig:Multilane}: Both second-order schemes provide a similar solution.

\section{Conclusion}\label{sec:4}
Based on the works \cite{kurganov2007reduction,kurganov2009non} we revisited two numerical schemes for nonlocal conservation laws: the KT scheme and the CU flux. For the latter one we provided a proof of convergence in the scalar case for first- and second-order schemes.
Although the KT scheme provides very accurate solutions, the computational costs are rather high.
In contrast, the CU flux is a very promising numerical flux function for nonlocal conservation laws. It achieves similar results as Godunov-type schemes (even coincides with them in many models) and for a second-order extension similar results as the KT scheme, while having less computational costs.
Since it is a central flux, it can be applied as a black-box solver, in particular, to very general systems of nonlocal conservation laws, for which Godunov-type schemes might be difficult to obtain.

\section*{Acknowledgment}
J.~F. is supported by the German Research Foundation (DFG) through SPP 2410 `Hyperbolic Balance Laws in Fluid Mechanics: Complexity, Scales, Randomness' under grant FR 4850/1-1. S.~R. is supported by NBHM, DAE, India (Ref. No. 02011/46/2021 NBHM(R.P.)/R \& D II/14874).
Part of this work was carried out during the research exchanges of J.F. and R.S. at the other institute as part of the PECFAR program. J.F. and R.S. would like to thank the IGSTC- DST, India and BMBF, Germany for their generous support.

\bibliographystyle{abbrv}
\bibliography{mybibiliography.bib}
\end{document}